\documentclass[preprint,12pt]{elsarticle}

\usepackage{graphicx}
\usepackage{color}
\usepackage{float}
\usepackage{amsfonts}
\usepackage{amssymb}
\usepackage{amsmath}
\usepackage{mathtools}
\usepackage{amssymb}
\usepackage{lineno}

\journal{http://www.arxiv.org}

\begin{document}

\begin{frontmatter}


\title{A Solution for the Open Abelian Sandpile Problem of Distributing k Items in N Vertices, where $k = N$}



\author{Michael Alexander Waddell}
\ead{mwadd002@fiu.edu}
\author{Edmundo Sebastian Barriga}
\ead{ebarr088@fiu.edu}

\address{Florida International University, United States}

\begin{abstract}
This paper outlines a closed solution to an open problem in Graph Theory concerning the classification of the successful initial distributions of k items in N vertices, where $\color{blue}k = N$, that lead to the terminal set $\color{blue}  N_k = \{ n_i\} $, where $\color{blue}n_i=1$ and $\color{blue}i=1,2,3,...,k$. First, each successful initial distribution is enumerated using an algorithm. The closed solution classifies the terminal set in terms of its modulus, and proves that each successful initial distribution can be classified by the same modulus.
\end{abstract}


\end{frontmatter}


\section{Introduction}
\label{S:1}

In the Abelian Sandpile Model, items are randomly assigned to discrete locations on the sandpile until a critical slope, or number of items. At this point, the locations that have a slope greater than the critical slope redistribute their items in a cascading fashion onto neighboring locations $\color{red} [1]$.   Consider an Abelian Sandpile described by a closed circle made of N vertices. Randomly distributed among the N vertices are k items, representing an initial distribution. With the goal of distributing the items among the vertices such that no vertex has more than 1 item placed on it, called the terminal set, any vertex A with more than 1 item has 2 items removed from it, with 1 item moving to the vertex right of vertex A and the other moving to the vertex to the left of vertex A. For $\color{blue} k < N $, the terminal set is reached after a finite number of turns and for $\color{blue} k > N $, the terminal set is never reached. For $\color{blue} k = N $, the initial distribution of the k items determines whether the terminal set is reached or not. 

\begin{figure}[H] 
\centering
\includegraphics[width=.75\textwidth]{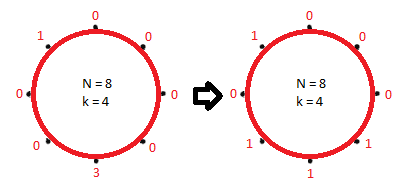}
\caption{k $<$ N Terminal set}
\end{figure}

\textbf{Open Problem:} Suppose $\color{blue} k = N$. Classify all successful initial distributions of items that lead to the terminal set. 

\section{Simplifying Cyclic Distributions}

Any circular set of vertices can be represented by a cyclic ordered set $\color{red} [2]$. For example, an $\color{blue} N = 3$ circle can be represented by the set: $\color{blue} N_3 = \{ n_1,n_2,n_3 \} $ and $\color{blue} \therefore N_k = \{ n_1,...,n_k \} $, where vertex $\color{blue}n_k$ is connected to vertex $\color{blue}n_1$ and vice-versa and $\color{blue}k$ is finite $\color{red} [3]$. Since each set is a cyclic order, the following is also true: 

$$\color{blue} N_3 = \{ n_1,n_2,n_3 \} = \{ n_3,n_1,n_2 \} = \{ n_2,n_3,n_1 \} ... $$

$$\color{blue} \therefore N_k = \{ n_1,...,n_k \} = \{ n_k,n_{k-1},..., n_1 \} = \{ n_{k-1},n_{k-2},...,n_1,n_k\} ... $$

, which is really just the equivalent of stating that the circle's order is independent of how you rotate it in space, and so the cyclic ordered set representing the circle. 

For any $\color{blue} k = N$, the terminal set is reached when an equal distribution of items among all of the vertices is achieved, i.e. where each vertex has 1 item.\footnote{ For instance, the terminal set $\color{blue} N_8 = \{ 1,1,1,1,1,1,1,1 \} $.} The first step in determining which initial distributions are successful is to enumerate all possible initial distributions between the set where all items are grouped in one vertex,  $\color{blue} N_k = \{ k,0,...,0 \} $, and the terminal set $\color{blue} N_k = \{ n_i\} $, where $\color{blue} n_i=1$ and $\color{blue} i=1,2,3,...,k$ .

Another way of representing the possible distributions of k  in N positions where 0's and permutations are allowed is the $\textbf{weak composition of k=N}$ $\color{red}[4]$. For example, all initial distributions for $\color{blue} k = N = 3 $ are:

$$\color{blue} W_3 = \{3,0,0\}, \{0,3,0\}, \{0,0,3\}, \{2,1,0\}, \{2,0,1\}, $$ 
$$\color{blue} \{ 1,2,0 \}, \{1,0,2\}, \{0,2,1\}, \{0,1,2\}, \{1,1,1\}  $$

The only initial distributions from this weak composition that reach the terminal set in a finite number of turns are: $\color{blue} \{3,0,0\}, \{0,3,0\}, \{0,0,3\}, \{1,1,1\}  $. These successful initial distributions yield the terminal set after a finite number of turns greater than or equal to 0 and thus include the terminal set. Using this system, the original algorithm is: $\color{blue} \{ (n_{i-1} + 1), (n_i - 2), (n_{i+1} + 1) \} $.

\section{Classifying the Enumeration of Successful Initial Distributions}

\textbf{Conjecture 1.0:} All successful initial distributions can be classified as the enumeration of all distributions that are a discrete number of reversed algorithmic turns away from the terminal set. A reversed algorithmic turn can be described as: $\color{blue} \{ (n_{i-1} - 1), (n_i + 2), (n_{i+1} - 1) \} $. 

 It is possible to see if each initial distribution reaches the terminal set in a finite number of turns with enough computing power by using the algorithm discussed in the introduction. However, it is impossible to know with the current method whether or not an arbitrary initial distribution will end in the terminal set.

 \textbf{Proof of Conjecture 1.0:} Every successful initial distribution has reached the terminal set in a finite $\color{blue}m_s$ turns, where $\color{blue}s$ is the particular successful starting position. It is true then that the non-terminal set described by $\color{blue} (m_s - 1)$  turns is a successful initial distribution that is 1 turn away from reaching the terminal set. This is also the case for the initial distribution that is $\color{blue} (m_s -1) - 1$  turns away from the terminal set.  In fact, since all turns are of the same type (any vertex with more than 1 item moves 1 to the left neighbor and 1 item to the right neighbor), by induction we can state that all successful initial distributions between the terminal set and the successful initial distribution $\color{blue} m_s$ can be described as being $\color{blue} (m_s - l_s)$ turns away from the terminal set, where $\color{blue} l_s \leq m_s$.

Assuming the successful initial distribution $\color{blue}m_s$ turns away from the terminal set allows for any more reversed algorithmic turns, the initial distribution described by $\color{blue} m_s + 1$, or the initial distribution that comes 1 turn before $\color{blue} m_s$ is also a valid successful initial distribution. This must also therefore be true for all successful initial distributions described by $\color{blue} (m_s + n_s)$ turns away from the terminal set, where $\color{blue} n_s \geq m_s$. In other words, all initial distributions $ m_s$ turns from the terminal set are successful initial distributions and therefore all possible successful initial distributions can be enumerated.

\textbf{Since this is true for any arbitrary successful initial distributions, all successful initial distributions can be classified in this way for the case where} $\color{blue} k = N$ \textbf{and no successful initial distribution is outside of this classification by definition.}

\section{Using the Classification to Derive a Closed Solution}

\textbf{Conjecture 2.0:} Every successful initial distribution shares the same characteristic modulus with the terminal set. 

\textbf{Proof of Conjecture 2.0:} For every 2-number item change in the number of items at vertex $\color{blue} n_i $, there has to be a 1-number item change in $\color{blue} n_{i-1} $  and $\color{blue}  n_{i+1} $. Therefore the initial distribution described by: $\color{blue} \{ n_{i-1}, n_i, n_{i+1} \} $, must be equivalent to the initial distribution: $\color{blue} \{ (n_{i-1} + 1), (n_i - 2), (n_{i+1} + 1) \} $, in that they are both successful initial distributions. It is not enough to know the change in magnitude of each vertex, but also the direction of the change in magnitude of the number of items at each vertex, and so classifying these as successful initial distributions must take both factors into account. Taking the sum of the product of the position of the vertex with the number of items on it, we can define the characteristic classification as: 

$$\color{blue}C_N = ( \sum_{i=1}^N i*n_i ) mod (N)  $$

, where $\color{blue} i = 1,2,3,...,N$ is the position in the initial distribution, and $\color{blue} n_i$ is the number of items at that vertex. This can be proven through basic algebra: 

$$\color{blue} C_N = (i_1*n_{i-1} + i_2*n_{i} + i_3*n_{i+1}) mod (N) $$

$$\color{blue} C_N' = [i_1*(n_{i-1}+1) + i_2*(n_{i}-2) + i_3*(n_{i+1}+1)] mod (N)  $$

$$\color{blue}  = (i_1n_{i-1} + i_1 + i_2n_i - 2i_2 + i_3n_{i+1} + i_3) mod (N)  $$

$$\color{blue}  = [(i_1n_{i-1} + i_2n_i + i_3n_{i+1})+( i_1 - 2i_2 + i_3)] mod (N)  $$

, where $\color{blue}  ( i_1 - 2i_2 + i_3) = (1-2+1) = 0  $ by definition since they are simply magnitude changes in position. Therefore:

$$\color{blue}  C_N = C_N' $$

Since this applies to any arbitrary $\color{blue}  n_i$ and therefore any $\color{blue}n_{i-1}, n_{i+1} $, and any change in the $\color{blue} m_s$ turns to acquire the terminal set must follow the $\color{blue} \{ (n_{i-1} + 1), (n_i - 2), (n_{i+1} + 1) \} $ algorithm in the cyclic ordered set, \textbf{this general proof applies to any successful initial distribution in relation to the terminal set and therefore any initial distribution can be checked to see if it ends in the terminal set.} Any initial distribution's classification that does not equal the terminal set's classification is therefore not a successful initial distribution. 

\newpage
\section{References}

[1]: Bak, P., Tang, C. and Wiesenfeld, K. (1987). "Self-organized criticality: an explanation of 1/ƒ noise". Physical Review Letters 59 (4): 381–384. Bibcode:1987PhRvL..59..381B. doi:10.1103/PhysRevLett.59.381.

[2]: Huntington, Edward V. (July 1935), "Inter-Relations Among the Four Principal Types of Order" (PDF), Transactions of the American Mathematical Society 38 (1): 1–9, doi:10.1090/S0002-9947-1935-1501800-1.

[3]: Courcelle, Bruno (21 August 2003), "2.3 Circular order", in Berwanger, Dietmar; Grädel, Erich, Problems in Finite Model Theory (PDF), p. 12.

[4]: Eger, Steffen (2013). "Restricted weighted integer compositions and extended binomial coefficients" . Journal of Integer Sequences 16.

\newpage
\section{Appendix}
In the Appendix, specific examples for each conjecture will be used to demonstrate their proofs. 

\subsection{Trivial Example of Conjecture 1.0}

For $\color{blue} k = N = 4$ the terminal set is $\color{blue} N_4 = \{ 1,1,1,1 \} $, and applying a reversed algorithmic turn to the 2nd vertex produces an initial distribution: 

$$\color{blue} N_4 = \{ (1-1),(1+2),(1-1),1 \} \rightarrow \{0,3,0,1 \} = m_s - 1 $$

, which is a valid successful initial distribution that ends in the terminal set in 1 step. Applying the same rule to the 3rd vertex which as a value of 0:

$$ \color{blue} N_4 = \{0,(3-1),(0+2),(1-1) \} \rightarrow \{0,2,2,0 \} = m_s - 2$$

, which is another valid successful initial distribution that ends in the terminal set in 2 steps, and so on.

\subsection{Trivial Example of Conjecture 2.0}

For $\color{blue} k = N = 10$. The terminal set is simply $\color{blue} N_{10} = \{ 1,1,1,1,1,1,1,1,1,1 \} $, and therefore: 

$$\color{blue} C_{10} = [1(1)+1(2)+1(3)+1(4)+1(5)+1(6)+1(7)+1(8)+1(9)+1(10)]mod(10) = 5  $$

A $\color{blue} m_s - 1$ successful initial distribution is $\color{blue} N_{10} = \{ 0,3,0,1,1,1,1,1,1,1 \} $, and the classification of this is: 

$$\color{blue} C_{10}' = [0(1)+3(2)+0(3)+1(4)+1(5)+1(6)+1(7)+1(8)+1(9)+1(10)]mod(10) = 5  $$

$$ \color{blue} \therefore C_{10} = C_{10}'  $$











\end{document}